\documentclass[11pt, a4paper]{amsart}

\usepackage{hyperref}
\usepackage[T1]{fontenc}
\usepackage{ae}

\renewcommand\AA{\mathbb{A}}
\newcommand\TT{\mathcal{T}}
\newcommand\OO{\mathcal{O}}
\newcommand\MM{\mathcal{M}}
\newcommand\NN{\mathcal{N}}
\newcommand\ZZ{\mathbb{Z}}
\newcommand\PP{\mathbb{P}}
\newcommand\LL{\mathcal{L}}
\newcommand\KK{\mathbb{K}}
\newcommand\RR{\mathbb{R}}
\newcommand\xx{\mathbf{x}}
\newcommand\Ptwo{{\PP^2}}
\newcommand\Pthree{{\PP^3}}
\newcommand{\Aone}{{\mathbf A}_1}
\newcommand{\Atwo}{{\mathbf A}_2}

\newcommand{\Afour}{{\mathbf A}_4}
\newcommand{\Esix}{{\mathbf E}_6}
\newcommand{\Eseven}{{\mathbf E}_7}
\newcommand{\Dfive}{{\mathbf D}_5}
\DeclareMathOperator{\NS}{NS}
\newcommand{\Tns}{{T_{\NS}}}
\DeclareMathOperator{\multiplikativ}{m}
\newcommand{\Gm}{\mathbb{G}_{\multiplikativ}}
\DeclareMathOperator{\Pic}{Pic}
\DeclareMathOperator{\Cox}{Cox}
\DeclareMathOperator{\Spec}{Spec}
\DeclareMathOperator{\Proj}{Proj}
\DeclareMathOperator{\tr}{tr}
\DeclareMathOperator{\rad}{rad}
\newcommand{\qa}[1]{q_{#1}}
\newcommand{\qb}[1]{q'_{#1}}
\newcommand{\qc}[1]{q''_{#1}}
\newcommand{\ep}{\epsilon}
\newcommand{\tH}{\widetilde H}
\newcommand{\Xiij}[2]{\Xi''_{#1,#2}}
\newtheorem{theorem}{Theorem}
\newtheorem{lemma}[theorem]{Lemma}
\newtheorem{proposition}[theorem]{Proposition}
\theoremstyle{definition}
\newtheorem{definition}[theorem]{Definition} 
\theoremstyle{remark}
\newtheorem{remark}[theorem]{Remark}

\begin{document}

\title[Universal torsors and homogeneous spaces]
{Universal torsors of Del Pezzo surfaces\\and homogeneous spaces}

\author{Ulrich Derenthal} 

\address{Mathematisches Institut, Universit\"at G\"ottingen,
  Bunsenstr. 3-5, 37073 G\"ottingen, Germany}

\email{derentha@math.uni-goettingen.de} 

\date{September 13, 2006}



\begin{abstract}
  Let $\Cox(S_r)$ be the homogeneous coordinate ring of the blow-up
  $S_r$ of $\Ptwo$ in $r$ general points, i.e., a smooth Del Pezzo
  surface of degree $9-r$.  We prove that for $r \in \{6,7\}$,
  $\Proj(\Cox(S_r))$ can be embedded into $G_r/P_r$, where $G_r$ is an
  algebraic group with root system given by the primitive Picard
  lattice of $S_r$ and $P_r\subset G_r$ is a certain maximal parabolic
  subgroup.
\end{abstract}

\maketitle

\tableofcontents

\section{Introduction}\label{sec:introduction}

In this note we continue our investigations of universal torsors over Del
Pezzo surfaces over an algebraically closed field $\KK$ of characteristic 0.
The blow-up $S_r$ of $\Ptwo$ in $r \le 8$ points in \emph{general
  position}\footnote{no three on a line, no six on a conic, no eight and one
  of them singular on a cubic curve} is a smooth Del Pezzo surface of degree
$9-r$; we will assume that $r \in \{3, \dots, 7\}$. A smooth Del Pezzo surface
of degree 3 (resp. degree 2) is a smooth cubic surface in $\Pthree$ (resp. a
double cover of $\Ptwo$ ramified in a smooth curve of degree 4). The Picard
group $\Pic(S_r)$ is a lattice with a non-degenerate symmetric linear form
$(\cdot,\cdot)$, the \emph{intersection form}. It is well-known that
$\Pic(S_r)$ contains a canonical root system $R_r$, which carries the action
of the associated Weyl group $W_r$, see Table~\ref{tab:weyl_groups} and
\cite{MR833513}.

\begin{table}[ht]
  \centering
  \[\begin{array}[h]{|c||c|c|c|c|c|}
    \hline
    r & 3 & 4 & 5 & 6 & 7 \\
    \hline\hline
    R_r & \Atwo+\Aone & \Afour & \Dfive & \Esix & 
    \Eseven\\
    \hline
    N_r & 6 & 10 & 16 & 27 & 56 \\
    \hline
  \end{array}\]
  \caption{The root systems associated to Del Pezzo surfaces.}
  \label{tab:weyl_groups}
\end{table}

It was a general expectation that the Weyl group symmetry on
$\Pic(S_r)$ should be a reflection of a geometric link between Del
Pezzo surfaces and algebraic groups.  Here we show that universal
torsors of smooth Del Pezzo surfaces of degree $2$ and $3$ admit an
embedding into a certain flag variety for the corresponding algebraic
group. The degree 5 case goes back to Salberger (talk at the
  Borel seminar Bern, June 1993) following Mumford \cite{MR0437531},
and independently Skorobogatov \cite{MR1260765}.  The degree 4 case
was treated in the thesis of Popov \cite[Chapter~6]{popov_diplom}.
The existence of such an embedding in general was conjectured by
Batyrev in his lecture at the conference \emph{Diophantine geometry}
(Universit\"at G\"ottingen, June 2004).  Skorobogatov announced
related work in progress (joint with Serganova) at the conference
\emph{Cohomological approaches to rational points} (MSRI, March 2006).

\medskip

As in \cite[Section~2]{MR2029863}, the simple roots of $R_{r-1}$ (with
$R_2 = \Atwo$) can be identified with a subset $I_r$ of the simple
roots of $R_r$ such that the edges in the Dynkin diagrams are
respected. The complement of $I_r$ in $R_r$ consists of exactly one
simple root $\alpha_r$, with associated fundamental weight $\varpi_r$.
Let $G_r$ be a simply connected linear algebraic group associated to
$R_r$, and fix a Borel group containing a maximal torus.  The Weyl
group $W_r$ acts on the weight lattice of $G_r$. The fundamental
representation $\varrho_r$ of $G_r$ with highest weight $\varpi_r$ has
dimension $N_r$ as listed in Table~\ref{tab:weyl_groups}; the weights
of $\varrho_r$ can be identified with classes of curves $E\subset S_r$
with self-intersection number $(E,E) = -1$, so called $(-1)$-curves.

Let $P_r$ be the maximal parabolic subgroup corresponding to $I_r$. By
\cite{MR541023}, we can regard the $G_r$-orbit $H_r$ of the weight space to
$\varpi_r$ as the affine cone over $G_r/P_r$, and $H_r$ is given by
quadratic equations in affine space $\AA^{N_r}$.  For $r=6$, the
equations are all partial derivatives of a certain cubic form on the
27-dimensional representation $\varrho_6$ of $G_6$.  For $r=7$,
equations can be found in \cite{MR0054609}. In both cases, the
equations were already known to E. Cartan in the 19th century. See
Section~\ref{sec:homogeneous_spaces} for more details on $G_r$,
$\varrho_r$, and $H_r$ for $r \in \{6,7\}$.

\medskip

The \emph{universal torsor} $\TT_r$ over $S_r$ is defined as follows:
Let $\LL_0, \dots, \LL_r$ be a basis of $\Pic(S_r) \cong \ZZ^{r+1}$,
and let $\LL_i^\circ := \LL_i \setminus \{\text{zero-section}\}$. Then
\[\TT_r := \LL_0^\circ \times_{S_r} \dots \times_{S_r} \LL_r^\circ.\]
It is a $\Tns(S_r)$-bundle over $S_r$, where $\Tns(S_r)$ is the
N\'eron-Severi torus of $S_r$.

The \emph{total coordinate ring}, or \emph{Cox ring} of $S_r$ is defined as
\[\Cox(S_r) := \bigoplus_{(\nu_0, \dots, \nu_r) \in \ZZ^{r+1}}
H^0(S_r,\LL_0^{\otimes \nu_0} \otimes \dots \otimes \LL_r^{\otimes
  \nu_r})\] as a vectorspace, and the multiplication is induced by the
multiplication of sections (see \cite{MR2001i:14059} and
\cite{MR2029863}). It is naturally graded by $\Pic(S_r)$, and it is
generated by $N_r$ sections corresponding to the $(-1)$-curves on
$S_r$. The ideal of relations in $\Cox(S_r)$ is generated by certain
quadratic relations which are homogeneous with respect to the
$\Pic(S_r)$-grading (see \cite{MR2029863} for more details). Let
$\AA(S_r) := \Spec(\Cox(S_r)) \subset \AA^{N_r}$ be the corresponding
affine variety. The universal torsor $\TT_r$ is an open subset of
$\AA(S_r)$ (cf.  \cite{MR2001i:14059}).  See \cite{MR2029868} and
\cite{math.AG/0604194} for the calculation of universal torsors and
Cox rings for \emph{singular} Del Pezzo surfaces and
\cite{math.AG/0603111} for the smooth cases.

\medskip

Universal torsors can be applied to \emph{Manin's conjecture} (see
\cite{MR89m:11060} and \cite{MR1032922}) on the number of rational
points of bounded height on Del Pezzo surfaces. For the moment, let
the Del Pezzo surface $S_r$ be defined over a number field $k$, and let
$H: S_r(k) \to \RR$ be the anticanonical height function. Let $U$ be the
complement of the $(-1)$-curves on $S_r$. Then Manin's conjecture
predicts that \[N_{U,H}(B) := \#\{\xx \in U(k) \mid H(\xx) \le B\}\]
behaves asymptotically as \[N_{U,H}(B) \sim c \cdot B \cdot (\log
B)^{n-1}\] for some positive constant $c$, where $n$ is the rank of
the Picard group (over $k$) of $S_r$.

In results concerning Manin's conjecture for various Del Pezzo
surfaces (see \cite{math.NT/0511041} and
\cite[Section~1]{math.NT/0604193} for an overview), often the first
step is a translation of the counting problem for rational points on
$S_r$ to the counting of integral points in certain ranges on a
universal torsor $\TT_r$. The number of these points on $\TT_r$ can
then be estimated using techniques from analytic number theory.

Salberger \cite{MR1679841} gave a proof of Manin's conjecture for
toric varieties, which include smooth Del Pezzo surfaces of degree
$\ge 6$, using universal torsors.  De la Bret\`eche \cite{MR1909606}
used Salberger's and Skorobogatov's description of the universal
torsor as a homogeneous space in his proof of the asymptotic formula
for a Del Pezzo surface of degree 5.  In lower degrees, Manin's
conjecture has been proved only for examples of singular Del Pezzo
surfaces (see \cite{math.NT/0509370} for a singular cubic surface,
using the universal torsor), but it is a general expectation that
universal torsors should lead to a proof of Manin's conjecture also in
the remaining smooth cases.

\medskip

We have seen that both $\AA(S_r)$ and $H_r$ can be viewed as embedded into
$\AA^{N_r}$, with a natural identification of the coordinates. For the
embedding of $\AA(S_r)$, we have some freedom: As the generators of
$\Cox(S_r)$ are canonical only up to a non-zero constant, we can choose a
rescaling factor for each of the $N_r$ coordinates, giving a $N_r$-parameter
family of embeddings of $\AA(S_r)$ into affine space.  The task is to find a
rescaling such that $\AA(S_r)$ is embedded into $H_r$.

More precisely, we start with an arbitrary embedding \[\AA(S_r)
\subset \AA_r := \Spec(\KK[\xi(E) \mid \text{$E$ is a $(-1)$-curve on
  $S_r$}]) \cong \AA^{N_r},\] and view \[H_r \subset \AA'_r :=
\Spec(\KK[\xi'(E) \mid \text{$E$ is a $(-1)$-curve on $S_r$}]) \cong
\AA^{N_r}\] as embedded into a different affine space. An isomorphism
$\phi_r : \AA_r \to \AA'_r$ such that \[\phi_r^*(\xi'(E)) =
\xi''(E)\cdot\xi(E)\] for each of the $N_r$ coordinates, with
$\xi''(E) \in \KK^*:=\KK \setminus \{0\}$, is called a
\emph{rescaling}, and the factors $\xi''(E)$ are (a system of)
\emph{rescaling factors}. A rescaling $\phi_r$ which embeds
$\AA(S_r)$ into $H_r$ is called a \emph{good rescaling}.

Our main result is:

\begin{theorem}\label{thm:main}
Let $S_r$ be a smooth Del Pezzo surface of degree $9-r$
and $\AA(S_r)$ the affine variety described above. 
Let $H_r$ be the affine
cone over the flag variety $G_r/P_r$ associated to the root system
$R_r$ as in Table~\ref{tab:weyl_groups}.

For $r \in \{6,7\}$, there exists a $(N_{r-2}+2)$-parameter family of
good rescalings $\phi_r$ which embed $\AA(S_r)$ into $H_r$.
\end{theorem}

\begin{remark}\label{rem:parameters}
  The number $N_{r-2}+2$ of parameters is 12 for $r=6$, respectively 18 for
  $r=7$.

  The rescaling factors are naturally graded by $\Pic(S_r) \cong \ZZ^{r+1}$,
  and we will see in Section~\ref{sec:rescalings} that the conditions for good
  rescaling are homogeneous with respect to this grading. Therefore, for each
  good rescaling $\phi_r$, there is a $(r+1)$-parameter family of good
  rescalings which differ from $\phi_r$ only by the action of $\Tns(S_r) \cong
  \Gm^{r+1}$.  Similarly, $\Tns(S_r)$ acts on $\AA(S_r)$, and it is easy to
  see that the image of $\AA(S_r)$ in $H_r$ is the same for all good
  rescalings in the same $(r+1)$-parameter family. Therefore, the
  $(N_{r-2}+2)$-parameter family of good rescalings gives rise to a
  $(N_{r-2}-r+1)$-parameter family (where $N_{r-2}-r+1$ equals 5 for $r=6$,
  resp.\ 10 for $r=7$) of images of $\AA(S_r)$ in $H_r$.

  For $r=5$, we have $N_{r-2}-r+1 = 2$, and by
  \cite[Section~6.3]{popov_diplom}, there is a two-parameter family of images
  of $\AA(S_5)$ under good rescalings in $H_5$.
\end{remark}

In Section~\ref{sec:del_pezzo}, we summarize results of \cite{MR2029863} and
\cite{math.AG/0603111} on Cox rings of Del Pezzo surfaces of degree 3 and 2. In
Section~\ref{sec:homogeneous_spaces}, we recall the classical equations for
the homogeneous spaces $G_r/P_r$ and give a simplified description on a
certain Zariski open subset; this will help to find good rescalings.  In
Section~\ref{sec:rescalings}, we derive conditions on the rescaling factors in
terms of the description of $\Cox(S_r)$ and $G_r/P_r$.  In
Section~\ref{sec:degree_3} and Section~\ref{sec:degree_2}, we determine good
rescalings in degree 3 and 2, finishing the proof of Theorem~\ref{thm:main}.

\medskip

\noindent\textbf{Acknowledgments.}  This work was completed at the
\emph{Mathematical Sciences Research Institute} (MSRI, Berkeley)
during the program \emph{Rational points on higher-dimensional
varieties} (Spring 2006). I am grateful for the invitation and ideal
working conditions. Calculations leading to these results were carried
out on computers of the \emph{Gauss-Labor} (Universit\"at G\"ottingen).

I am grateful to my advisor Yu.\ Tschinkel for introducing me to this
problem.  I have benefitted from conversations and correspondence with
V.\ Batyrev and B.\ Hassett. I thank the anonymous referee for several
helpful suggestions.

\section{Cox rings of Del Pezzo surfaces}\label{sec:del_pezzo}

In this section, we describe the Cox ring of Del Pezzo surfaces of
degree 3 and 2 in order to fix some notation. The results can be found
in \cite{MR2029863} and \cite{math.AG/0603111}. 

Let $r \in \{6,7\}$. Without loss of generality, we may assume that four of
the $r$ blown-up points in $\Ptwo$ giving $S_r$ are in the positions
\begin{equation*}
  \begin{split}
    p_1 &= (1:0:0), \qquad p_2 = (0:1:0), \qquad p_3 = (0:0:1),\qquad
    p_4 = (1:1:1).
  \end{split}
\end{equation*}
By \cite[Theorem~3.2]{MR2029863}, the generators of $\Cox(S_r)$ are
sections $\xi(E)$ vanishing in a $(-1)$-curve $E$ on $S_r$. Their
number is $N_r$ as in Table~\ref{tab:weyl_groups}. We use the same
symbols $\xi(E)$ for the coordinates in $\AA_r$. Let $-K_r$ be the
anticanonical divisor class of $S_r$.

A $(k)$-ruling $D \in \Pic(S_r)$ is the sum of two $(-1)$-curves whose
intersection number is $k$, cf.  \cite[Definition~4.6]{MR2029863} and
\cite[Definition~1]{math.AG/0603111}. The relations come in groups of
$r-3$ for each $(1)$-ruling $D$. We will denote them by \[F_{D,1},
\dots, F_{D,r-3}.\] For $r=7$, we have further relations corresponding
to the $(2)$-ruling $-K_7$.  The ideal $J_r$ generated by these
relations defines $\AA(S_r)$ (see \cite[Theorem~4.9]{MR2029863} for $r
\le 6$ and \cite[Theorem~2]{math.AG/0603111} for $r=7$).

\medskip

For $r = 6$, we assume \[p_5 = (1:a:b),\qquad p_6 = (1:c:d).\] The
$N_6=27$ $(-1)$-curves $E$ are denoted by $E_i$, $m_{i,j}$, and $Q_i$
corresponding to the six blown-up points, the 15 transforms of the
lines through two of the six points, and the six transforms of the
conics through five points as described in \cite[Section~3]{math.AG/0603111}. Let $\eta_i := \xi(E_i)$,
$\mu_{i,j}:=\xi(m_{i,j})$, and $\lambda_i:=\xi(Q_i)$. We order them in
the following way:
\[\eta_1, \dots, \eta_6, \quad \mu_{1,2}, \dots, \mu_{1,6}, \mu_{2,3},
\dots, \mu_{2,6}, \mu_{3,4}, \dots, \mu_{5,6}, \quad \lambda_1, \dots,
\lambda_6.\] The $(1)$-rulings are $-K_6-E$, where $E$ runs through
the $(-1)$-curves. We have
\[F_{-K_6-E,1}=\qa E, \qquad F_{-K_6-E,2}=\qb E, \qquad
F_{-K_6-E,3}=\qc E,\] where the 81 equations $\qa E, \qb E, \qc E$ are
listed in \cite[Section~3]{math.AG/0603111}.

\medskip

For $r=7$, let \[p_5 = (1:a_1:b_1), \qquad p_6 = (1:a_2:b_2), \qquad p_7 =
(1:a_3:b_3).\] The $N_7=56$ $(-1)$-curves $E_i, m_{i,j}, Q_{i,j}, C_i$ and the
corresponding generators $\xi(E)$ are described in \cite[Section~4]{math.AG/0603111}. They are ordered as
\[\eta_1, \dots, \eta_7, \mu_{1,2}, \dots, \mu_{1,7}, \mu_{2,3},
\dots, \mu_{6,7}, \nu_{1,2}, \dots, \nu_{1,7}, \nu_{2,3}, \dots,
\nu_{6,7}, \lambda_1, \dots, \lambda_7.\] By \cite[Theorem~2]{math.AG/0603111}, there are 529 relations $F_{D,i}$ (four for each
of the 126 $(1)$-rulings, and 25 for the $(2)$-ruling $-K_7$) between
them.  We do not want to list them here, as they can be determined by
the method of \cite[Lemma~4]{math.AG/0603111} in a straightforward
manner.  

\section{Homogeneous spaces}\label{sec:homogeneous_spaces}

In this section, we examine the equations defining the affine cone
$H_r \subset \AA'_r$ over $G_r/P_r$ for $r \in \{6,7\}$. For the $N_r$
coordinates $\xi'(E)$ of $\AA'_r$, we also use the names $\eta'_i$,
$\mu'_{i,j}$, $\lambda'_i$, and furthermore $\nu'_{i,j}$ in the case
$r=7$, with the obvious correspondence to the coordinates of $\AA_r$
as in the previous section.

In particular, we show that $H_r$ is a complete intersection on the
open subset $U_r$ of $\AA'_r$ where the coordinates $\eta'_1, \dots,
\eta'_r$ are non-zero.

\medskip

We will see that $H_r$ is defined by quadratic relations which are
homogeneous with respect to the $\Pic(S_r)$-grading. For each
$(1)$-ruling $D$, we have exactly one relation $p_D$ of degree $D$,
and furthermore in the case $r=7$, we have eight relations
$p^{(1)}_{-K_7}, \dots, p^{(8)}_{-K_7}$ where we use the convention
$p_{-K_7} := p^{(1)}_{-K_7}$. For any possibility to write $D$ as the
sum of two $(-1)$-curves $E$, $E'$, the relation $p_D$ has a term
$\xi'(E)\xi'(E')$ with a non-zero coefficient.

\begin{definition}\label{def:subsets}
  For a $(-1)$-curve $E$, let $U_E$ be the open subset of $\AA'_r$
  where $\xi'(E)$ is non-zero. Let $\NN(E)_k$ be the set of
  $(-1)$-curves $E'$ with $(E,E') = k$, and let $\Xi'(E)_k$ be the set
  of the corresponding $\xi'(E')$. Let $\NN(E)_{>k}$ and
  $\Xi'(E)_{>k}$ be defined similarly, but with the condition $(E,E')
  > k$.

  Let $U_r \subset \AA'_r$ be the intersection of $U_{E_1}, \dots,
  U_{E_r}$.
\end{definition}

Note that $\NN(E)_0$ has exactly $N_{r-1}$ elements because we can
identify its elements with the $(-1)$-curves on $S_{r-1}$. Since the
only $(-1)$-curve intersecting $E$ negatively is $E$ itself, the
number of elements of $\NN(E)_{>0}$ is $N_r-N_{r-1}-1$.

\begin{proposition}\label{prop:homspace}
  Let \[\Phi_r: H_r \cap U_r \to U_{r-1} \times (\AA^1 \setminus
  \{0\})\] be the projection to the coordinates $\xi'(E) \in
  \Xi'(E_1)_0$ and $\eta'_1$. The map $\Phi_r$ is an isomorphism. The
  dimension of $H_r$ is $N_{r-1}+1$.
\end{proposition}

\begin{proof}
  If $D=E_1+E$ is a $(1)$-ruling, then $(E_1,D)=0$, and all variables
  occurring in $p_D$ besides $\eta'_1$ and $\xi'(E)$ are elements of
  $\Xi'(E_1)_0$. For $\eta'_1 \ne 0$, the relation $p_D$ expresses
  $\xi'(E)$ in terms of $\eta'_1$ and $\Xi'(E_1)_0$.

  For a $(2)$-ruling $D=E_1+E$, we have $(E_1,D)=1$, so the relation
  $p_D$ expresses $\xi'(E)$ in terms of $\eta'_1$ and monomials
  $\xi'(E_i')\xi'(E_i'')$ where $\xi(E_i') \in \Xi'(E_1)_0$ and
  $\xi'(E_i'') \in \Xi'(E_1)_1$. Using the expressions for the
  elements of $\Xi'(E_1)_1$ of the first step, this shows that we can
  express the coordinates $\Xi'(E_1)_{>0}$ in terms of $\eta'_1$ and
  $\Xi'(E_1)_0$ by using the $N_r-N_{r-1}-1$ relations $g_{E_1+E}$ for
  $E \in \NN(E)_{>0}$.  This allows us to construct a map \[\Psi_r:
  U_{r-1} \times (\AA^1 \setminus \{0\}) \to \AA'_r.\]

  It remains to show that the image of $\Psi_r$ is in $H_r$, i.e.,
  that the resulting point also satisfies the remaining equations
  which define $H_r$. This is done in Lemma~\ref{lem:homspace_degree3}
  and Lemma~\ref{lem:homspace_degree2} below.
\end{proof}

\begin{remark}
  Proposition~\ref{prop:homspace} is also true if we enlarge $U_r$ to
  $U_{E_1}$ and $U_{r-1}$ to $\AA'_{r-1}$. However, the proofs of
  Lemma~\ref{lem:homspace_degree3} and
  Lemma~\ref{lem:homspace_degree2} are slightly simplified by
  restricting to $U_r$.
\end{remark}

First, we consider the case $r=6$. Consider the cubic form in $N_6=27$
variables \[F(M_1, M_2, M_3) := \det M_1+\det M_2+\det M_3-
\tr(M_1 M_2 M_3),\] where
\[M_1 :=
\begin{pmatrix}
  \eta'_1 & \lambda'_1 & \mu'_{2,3}\\
  \eta'_2 & \lambda'_2 & \mu'_{1,3}\\
  \eta'_3 & \lambda'_3 & \mu'_{1,2}\\
\end{pmatrix}, \qquad M_2:=
\begin{pmatrix}
  \lambda'_4 & \lambda'_5 & \lambda'_6\\
  \eta'_4 & \eta'_5 & \eta'_6\\
  \mu'_{5,6} & \mu'_{4,6} & \mu'_{4,5}\\
\end{pmatrix},\] and \[M_3 :=
\begin{pmatrix}
  \mu'_{1,4} & \mu'_{2,4} & \mu'_{3,4}\\
  \mu'_{1,5} & \mu'_{2,5} & \mu'_{3,5}\\
  \mu'_{1,6} & \mu'_{2,6} & \mu'_{3,6}\\
\end{pmatrix}.\] By \cite[Proposition~1.6]{MR2235344}, the group of invertible
$N_6 \times N_6$-matrices which leave invariant $F$ is a simply connected
linear algebraic group $G_6$ of type $\Esix$.

Note that the terms of $\tr(M_1M_2M_3)$ are
$M_1^{(i,j)}M_2^{(j,k)}M_3^{(k,i)}$ for $i,j,k \in \{1,2,3\}$ (where
$M_a^{(b,c)}$ is the entry $(b,c)$ of the matrix $M_a$), so the number of
terms of $F$ is $3\cdot 6+3^3=45$. Each is a product of three variables
$\xi'(E)$, $\xi'(E')$, $\xi'(E'')$ such that the corresponding $(-1)$-curves
$E$, $E'$, $E''$ on $S_6$ form a triangle, and their divisor classes add up to
$-K_6$. The coefficient is $+1$ in the nine cases
\begin{equation*}
  \begin{array}[h]{ccc}
    \eta'_1\mu'_{1,2}\lambda'_2, &\eta'_2\mu'_{2,3}\lambda'_3,
    &\eta'_3\mu'_{1,3}\lambda'_1,\\ \eta'_4\mu'_{4,6}\lambda'_6,
    &\eta'_5\mu'_{4,5}\lambda'_4, &\eta'_6\mu'_{5,6}\lambda'_5,\\
    \mu'_{1,4}\mu'_{2,5}\mu'_{3,6}, &\mu'_{1,5}\mu'_{2,6}\mu'_{3,4},
    &\mu'_{1,6}\mu'_{2,4}\mu'_{3,5}
  \end{array}
\end{equation*}
and $-1$ in the remaining 36 cases. (Of course, there is some choice
here, for example by permuting the indices $1, \dots, 6$, but it is
not as simple as choosing any 9 of the 45 terms to have the coefficient
$+1$. See \cite[Section~5]{MR1854697} for more details.)  

Let $\alpha_6$ be the simple root at the end of one of the ``long
legs'' in the Dynkin diagram $\Esix$. Let $\varpi_6$ be the associated
fundamental weight. The action of $G_6$ on $\KK^{N_6}$ is a
$N_6$-dimensional irreducible representation of $G_6$ whose highest
weight is $\varpi_6$ (cf. \cite[Section~20.2]{MR0106966}). The orbit
$H_6$ of the weight space of $\varpi_6$ is described by the vanishing
of the $N_6$ partial derivatives of the cubic form $F$ (see
\cite[Section~III.2.5]{MR1234494}).

The derivative with respect to $\xi'(E)$ contains five terms $\pm
\xi'(E')\xi'(E'')$ corresponding to the five ways to write the $(1)$-ruling
$D:=-K_6-E$ as the sum of two intersecting $(-1)$-curves $E', E''$. We will
denote it by $p_D = p_{-K_6-E}$.

\begin{lemma}\label{lem:homspace_degree3}
  For $\eta'_1 \ne 0$ and any values of \[\Xi'(E_1)_0 = \{\eta'_2,
  \dots, \eta'_6, \mu'_{2,3}, \dots, \mu'_{5,6}, \lambda'_1\}\] with
  non-zero $\eta'_2, \dots, \eta'_6$, the equations $p_{E_1+E}$ for
  \[E \in \NN(E_1)_1 = \{m_{1,2}, \dots, m_{1,6}, Q_2, \dots, Q_6\}\]
  define a point of $H_6$.
\end{lemma}

\begin{proof}
  As $\Tns(S_6)$ acts on $H_6$ and $\{E_1, \dots, E_6\}$ is a subset of a
  basis of $\Pic(S_6)$, we may assume that $\eta'_1 = \dots = \eta'_6 = 1$.
  Then for $i \in \{2, \dots, 6\}$, the equation $p_{E_1+m_{1,i}}$ allows us to
  express $\mu'_{1,i}$ in terms of the remaining $\mu'_{i,j}$:
  \begin{equation*}
    \begin{split}
    \mu'_{1,2}&= \mu'_{2,3} + \mu'_{2,4} + \mu'_{2,5} + \mu'_{2,6},\qquad
    \mu'_{1,3}= \mu'_{2,3} - \mu'_{3,4} - \mu'_{3,5} - \mu'_{3,6},\\
    \mu'_{1,4}&= -\mu'_{2,4} - \mu'_{3,4} + \mu'_{4,5} - \mu'_{4,6},\qquad
    \mu'_{1,5}= -\mu'_{2,5} - \mu'_{3,5} - \mu'_{4,5} + \mu'_{5,6}\\
    \mu'_{1,6}&= -\mu'_{2,6} - \mu'_{3,6} + \mu'_{4,6} - \mu'_{5,6}
    \end{split}
  \end{equation*}
  Furthermore, for $i \in \{2, \dots, 6\}$, we can use $p_{E_1+Q_i}$ in
  order to express $\lambda'_i$ in terms of $\lambda'_1$ and $\mu'_{j,k}$:
  \begin{equation*}
    \begin{split}
    \lambda'_2&=\mu'_{3,4}\mu'_{5,6} + \mu'_{3,5}\mu'_{4,6} + \mu'_{3,6}\mu'_{4,5} + \lambda'_1\\
    \lambda'_3&=-\mu'_{2,4}\mu'_{5,6} - \mu'_{2,5}\mu'_{4,6} - \mu'_{2,6}\mu'_{4,5} + \lambda'_1\\
    \lambda'_4&=-\mu'_{2,3}\mu'_{5,6} + \mu'_{2,5}\mu'_{3,6} - \mu'_{2,6}\mu'_{3,5} - \lambda'_1\\
    \lambda'_5&=-\mu'_{2,3}\mu'_{4,6} - \mu'_{2,4}\mu'_{3,6} + \mu'_{2,6}\mu'_{3,4} - \lambda'_1\\
    \lambda'_6&=-\mu'_{2,3}\mu'_{4,5} + \mu'_{2,4}\mu'_{3,5} - \mu'_{2,5}\mu'_{3,4} - \lambda'_1
    \end{split}
  \end{equation*}
  By substituting and expanding, we check that the remaining 17 relations
  are fulfilled. Therefore, the resulting point lies in $H_6$.
\end{proof}

\medskip

Next, we obtain similar results in the case $r=7$ with $N_7=56$. By
\cite[Corollary~2.6]{MR2235344}, a simply connected linear algebraic
group $G_7$ of type $\Eseven$ is obtained as the identity component of
the group of invertible $N_7 \times N_7$-matrices which leave
invariant a certain quartic form defined on a vectorspace of dimension
$N_7$ as in \cite[Section~2.1]{MR2235344}. The action of $G_7$ on this
vectorspace is an irreducible representation whose highest weight
$\varpi_7$ is the fundamental weight corresponding to the simple root
$\alpha_7$ at the end of the ``longest leg'' of the Dynkin diagram
$\Eseven$ (cf. \cite[Section~20.2]{MR0106966}).

We describe the orbit $H_7$ of the weight space of $\varpi_7$ under
$G_7$ explicitly. The $N_7$ coordinates $\xi'(E)$ in $\AA'_7$ are $\eta_i',
\mu_{j,k}', \nu_{j,k}', \lambda'_i$ for $i,j,k \in \{1, \dots, 7\}$ and $j <
k$.  The equations for $H_7$ are described in \cite{MR0054609} in terms of 56
coordinates $x^{ij}, y_{ij}$ ($i<j\in \{1, \dots, 8\}$). They correspond to
our variables as follows:
\[\eta_i' = x^{i8}, \qquad \mu_{k,l}' = y_{kl}, \qquad \nu_{k,l}' =
x^{kl}, \qquad \lambda_i' = y_{i8}.\]

For the $(1)$-rulings $D$ \cite[Lemma~11]{math.AG/0603111}, the relations
$p_D$ are $u^{ijkl}$ and $v^i_j$ as below. In the first column of
Table~\ref{tab:relations_degree2}, we list a symbol $D^{(n)}_I$ assigned to
the $(1)$-ruling in the second column, and the third column gives
the corresponding relation.

\begin{table}[ht]
  \centering
  \[\begin{array}[h]{|c||c|c|}
    \hline
    \text{symbol} & \text{$(1)$-ruling $D=D^{(n)}_I$} & \text{relation $p_D$}\\
    \hline\hline
    D^{(1)}_i & H-E_i &  v^8_i\\
    D^{(2)}_{i,j,k} & 2H-(E_1+\dots+E_7)+E_i+E_j+E_k & u^{ijk8}\\
    D^{(3)}_{i,j} & 3H-(E_1+\dots+E_7)+E_i-E_j & v^i_j\\
    D^{(4)}_{i,j,k,l} & 4H-2(E_1+\dots+E_7)+E_i+E_j+E_k+E_l & u^{ijkl}\\
    D^{(5)}_i & 5H-2(E_1+\dots+E_7)+E_i & v^i_8\\
    \hline
  \end{array}\]
  \caption{Rulings and relations defining $G_7/P_7$.}
  \label{tab:relations_degree2}
\end{table}

Let
\[u^{ijkl} = x^{ij}x^{kl}-x^{ik}x^{jl}+x^{il}x^{jk}+
\sigma\cdot(y_{ab}y_{cd}-y_{ac}y_{bd}+y_{ad}y_{bc}),\] where
$i<j<k<l$ and $a<b<c<d$, with $(i,j,k,l,a,b,c,d)$ a permutation of
$(1,\dots,8)$, and $\sigma$ its sign. For $i \ne j$
\[v^i_j = \sum_{k \in (\{1, \dots, 8\} \setminus \{i,j\})} x^{ik}y_{kj},\] 
where $x^{ba}=-x^{ab}$ and $y_{ba}=-y_{ab}$ if $b>a$.

For the $(2)$-ruling $-K_7$, we have the following eight equations
with 28 terms:
\[p^{(i)}_{-K_7} := v^i_i := -\frac 3 4 \sum_{j \in (\{1, \dots, 8\} \setminus \{i\})}
x^{ij}y_{ij} + \frac 1 4 \sum_{j<k \in (\{1, \dots, 8\} \setminus
  \{i\})} x^{jk}y_{jk}\]

\begin{lemma}\label{lem:homspace_degree2}
  For $\eta'_1, \dots, \eta'_7 \ne 0$, the 28 coordinates
  \[\eta'_i\quad(i \in \{1, \dots, 7\}),\quad \mu'_{j,k}\quad (j < k \in \{2,
  \dots, 7\}),\quad \nu'_{1,l} \quad (l \in \{2, \dots, 7\})\] in
  $\Xi'(E_1)_0$ and the 28 equations $p_D$ for \[D \in \{D^{(1)}_2,
  \dots, D^{(1)}_7, D^{(2)}_{1,2,3}, \dots, D^{(2)}_{1,6,7},
  D^{(3)}_{1,2}, \dots, D^{(3)}_{1,7}, -K_7\}\] define
  \[\mu'_{1,i} \quad (i \in \{2, \dots, 7\}),\quad \nu'_{j,k} \quad
  (j<k \in \{2, \dots, 7\}),\quad \lambda'_l\quad (l \in \{1, \dots,
  7\}),\] resulting in a point on $H_7$.

  Furthermore, we may replace $p_{-K}$ by $p_D$ for $D=D^{(3)}_{2,1}$.
\end{lemma}

\begin{proof}
  As above, we may assume that $\eta'_1= \dots= \eta'_7=1$ because of
  the action of $\Tns(S_7)$. For the 27 $(-1)$-curves $E \in
  \NN(E_1)_1$, the equation $p_{E_1+E}$ defines $\xi'(E)$ directly in
  terms of the 28 variables in $\Xi'(E_1)_0$; we do not list the
  expressions here. By substituting these results, we use $v^1_1$ in
  order to express $\lambda'_1$ in terms of these variables:
  \begin{equation*}
    \begin{split}
      \lambda'_1=&-\mu'_{2,3}\mu'_{4,5}\mu'_{6,7} + \mu'_{2,3}\mu'_{4,6}\mu'_{5,7} - \mu'_{2,3}\mu'_{4,7}\mu'_{5,6} + \mu'_{2,4}\mu'_{3,5}\mu'_{6,7} -
      \mu'_{2,4}\mu'_{3,6}\mu'_{5,7}\\ &+ \mu'_{2,4}\mu'_{3,7}\mu'_{5,6} - \mu'_{2,5}\mu'_{3,4}\mu'_{6,7} + \mu'_{2,5}\mu'_{3,6}\mu'_{4,7} - \mu'_{2,5}\mu'_{3,7}\mu'_{4,6} +
      \mu'_{2,6}\mu'_{3,4}\mu'_{5,7}\\ &- \mu'_{2,6}\mu'_{3,5}\mu'_{4,7} + \mu'_{2,6}\mu'_{3,7}\mu'_{4,5} - \mu'_{2,7}\mu'_{3,4}\mu'_{5,6} + \mu'_{2,7}\mu'_{3,5}\mu'_{4,6} -
      \mu'_{2,7}\mu'_{3,6}\mu'_{4,5}\\
      &- \mu'_{2,3}\lambda'_2 + \mu'_{2,3}\lambda'_3 - \mu'_{2,4}\lambda'_2 + \mu'_{2,4}\lambda'_4 - \mu'_{2,5}\lambda'_2 + \mu'_{2,5}\lambda'_5\\& - \mu'_{2,6}\lambda'_2 +
      \mu'_{2,6}\lambda'_6 - \mu'_{2,7}\lambda'_2 + \mu'_{2,7}\lambda'_7 - \mu'_{3,4}\lambda'_3 + \mu'_{3,4}\lambda'_4\\ &- \mu'_{3,5}\lambda'_3 + \mu'_{3,5}\lambda'_5 -
      \mu'_{3,6}\lambda'_3 + \mu'_{3,6}\lambda'_6 - \mu'_{3,7}\lambda'_3 + \mu'_{3,7}\lambda'_7\\ &- \mu'_{4,5}\lambda'_4 + \mu'_{4,5}\lambda'_5 - \mu'_{4,6}\lambda'_4 +
      \mu'_{4,6}\lambda'_6 - \mu'_{4,7}\lambda'_4 + \mu'_{4,7}\lambda'_7\\ &- \mu'_{5,6}\lambda'_5 + \mu'_{5,6}\lambda'_6 - \mu'_{5,7}\lambda'_5 + \mu'_{5,7}\lambda'_7 -
      \mu'_{6,7}\lambda'_6 + \mu'_{6,7}\lambda'_7
    \end{split}
  \end{equation*}
  We check directly by substituting and expanding that the remaining
  equations defining $H_7$ are fulfilled.

  As $v^2_1$ contains the term $\eta'_2\lambda'_1$, and $\eta'_2 \ne
  0$, we may replace $v^1_1$ by $v^2_1$.
\end{proof}

\section{Rescalings}\label{sec:rescalings}

Let $r \in \{6,7\}$. We follow the strategy of the case $r=5$
\cite[Section~6.3]{popov_diplom} in order to describe conditions for good
rescalings explicitly in terms of the rescaling factors. However, we use the
results of the previous section to simplify this
as follows:

Let \[\MM_6:=\{E_1+E \mid E \in \NN(E_1)_1\}\] and let
\[\MM_7:=\{E_1+E \mid E \in \NN(E_1)_1\} \cup \{D^{(3)}_{2,1}\}.\]
Let $\tH_r \subset \AA'_r$ be the variety defined by the equations
$g_D$ for $D \in \MM_r$.

By Proposition~\ref{prop:homspace}, Lemma~\ref{lem:homspace_degree3},
and Lemma~\ref{lem:homspace_degree2}, $H_r \cap U_r = \tH_r \cap U_r$.

\begin{remark}
  Because of $\NN(E_1)_2 = \{C_1\}$, it could be considered more
  natural to use $-K_7=E_1+C_1$ instead of $D^{(3)}_{2,1}=E_2+C_1$ in
  the definition of $\MM_7$. However, we choose to avoid the
  $(2)$-ruling $-K_7$ for technical reasons.
\end{remark}

\begin{lemma}
  A rescaling $\phi_r:\AA_r \to \AA'_r$ is good if and only if it
  embeds $\AA(S_r)$ into $\tH_r$.
\end{lemma}

\begin{proof}
  As $H_r \subset \tH_r$, a good rescaling $\phi_r$ satisfies
  $\phi_r(\AA(S_r)) \subset \tH_r$. Conversely, we have
  \[\phi_r(\AA(S_r)) \cap U_r \subset \tH_r \cap U_r = H_r \cap U_r\]
  by Lemma~\ref{lem:homspace_degree3} and
  Lemma~\ref{lem:homspace_degree2}.  Taking the closure and using that
  $H_r$ is closed and that $\AA(S_r)$ is irreducible by
  \cite{MR2029863}, we conclude that $\phi_r(\AA(S_r)) \subset H_r$,
  so the rescaling is good.
\end{proof}

As in Section~\ref{sec:del_pezzo}, let $J_r$ be the ideal
defining $\AA(S_r)$ in $\AA_r$.

In terms of the coordinate rings $\KK[\AA_r]$ and $\KK[\AA'_r]$ and in
view of the previous lemma, a rescaling $\phi_r$ is good if, for all $D \in
\MM_r$, the ideal $J_r \subset \rad(J_r)$ contains $\phi_r^*(p_D)$,
where $p_D$ is the equation defining $H_r$ corresponding to the
$(1)$-ruling $D$.

As $\KK[\AA_r]$ and $\KK[\AA'_r]$ are both graded by $\Pic(S_r)$ and
$\phi_r^*$ respects this grading, we need rescaling
factors such that $\phi_r^*(p_D)$ of degree $D \in \MM_r$ is a linear
combination of the equations $F_{D,1}, \dots, F_{D,r-3} \in J_r$. For
concrete calculations in the next sections, we describe this more
explicitly:

Let $D \in \MM_r$ be a $(1)$-ruling, which can be written in $r-1$
ways as the sum of two $(-1)$-curves $E_i', E_i''$. For $i \in \{1,
\dots, r-1\}$, let
\[\xi_i:=\xi(E_i')\xi(E_i''), \qquad \xi_i':=\xi'(E_i')\xi'(E_i''), \qquad
\xi_i'':=\xi''(E_i')\xi''(E_i'').\] Then $p_D$ has the form
\begin{equation}\label{eq:defining_Hr}
  p_D = \sum_{i=1}^{r-1}\ep_i\xi'_i
\end{equation}
with $\ep_i \in \{\pm 1\}$.

As  $\xi_i$ vanishes exactly on $E_i' \cup E_i''$, 
the 2-dimensional space $H^0(S_r, \OO(D))$ is generated
by any two $\xi_i$, $\xi_{i'}$. Hence, all other $r-3$ elements
$\xi_j$ are linear combinations of $\xi_i$, $\xi_{i'}$, with
non-vanishing coefficients. This gives $r-3$ relations of degree $D$
in $\Cox(S_r)$. Rearranging $\xi_1, \dots, \xi_{r-1}$ such that the
two elements $\xi_i$, $\xi_{i'}$ of our choice have the indices $r-2$
and $r-1$, we can write them as
\begin{equation}\label{eq:defining_FDj}
  F_{D,j} = \xi_j + \alpha_j \xi_{r-2} + \beta_j \xi_{r-1},
\end{equation}
for $j \in \{1, \dots, r-3\}$, where $\alpha_j, \beta_j \in \KK^*$.

Suppose that $\phi_r^*(p_D)$ is a linear combination of the $F_{D,j}$
with factors $\lambda_j$: \[\phi_r^*(p_D) - \sum_{j=1}^{r-3} \lambda_j
F_{D,j}=0.\] Since $\phi_r^*(\xi'(E)) = \xi''(E)\cdot\xi(E)$, we have
$\phi_r^*(\xi_i') = \xi_i''\cdot\xi_i$ for the monomials of degree 2.
Then the above equation is equivalent to the vanishing of
\[\sum_{i=1}^{r-3} (\ep_i\xi_i'' - 
\lambda_i)\xi_i + \left(\ep_{r-2}\xi_{r-2}''-
  \sum_{j=1}^{r-3}\lambda_j\alpha_j\right)\xi_{r-2} +
\left(\ep_{r-1}\xi_{r-1}''-
  \sum_{j=1}^{r-3}\lambda_j\beta_j\right)\xi_{r-1}.\]

For $i \in \{1, \dots, r-3\}$, we see by considering the coefficients
of $\xi_i$ that we must choose $\lambda_i = \ep_i \xi_i''$. With this,
consideration of the coefficients of $\xi_{r-2}$ and $\xi_{r-1}$
results in the following conditions $g_{D,1}$, $g_{D,2}$ on the
rescaling factors $\xi_j''$, which are homogeneous of degree $D \in
\Pic(S_r)$:
\[ g_{D,1} := \ep_{r-2} \xi_{r-2}''-\sum_{j=1}^{r-3} \ep_j\alpha_j
\xi_j'' = 0,\qquad g_{D,2} := \ep_{r-1} \xi_{r-1}''-\sum_{j=1}^{r-3}
\ep_j\beta_j \xi_j'' = 0.\] Note that our choice of $\xi_{r-2}$ and
$\xi_{r-1}$ in the definition of $F_{D,j}$ as discussed before
\eqref{eq:defining_FDj} is reflected here in the sense that $g_{D,1}$
and $g_{D,2}$ express the corresponding $\xi''_{r-2}$ and
$\xi''_{r-1}$ as linear combinations of $\xi''_1, \dots, \xi''_{r-3}$
with non-zero coefficients.

This information can be summarized as
follows:

\begin{lemma}\label{lem:equations}
  For $r \in \{6,7\}$, a rescaling is good if and only if the rescaling
  factors $\xi''(E)$ fulfill the equations $g_{D,1}$ and $g_{D,2}$ for each
  $(1)$-ruling $D \in \MM_r$.
  
  As described above precisely, the non-zero coefficients $\ep_i$ are taken
  from the equations $p_D$ \eqref{eq:defining_Hr} defining $H_r$, and the
  non-zero $\alpha_j$, $\beta_j$ are taken from the equations $F_{D,j}$
  \eqref{eq:defining_FDj} defining $\AA(S_r)$.
\end{lemma}

Let $\Xi''(E)_k$ (resp. $\Xi''(E)_{>k}$) be the set of all
$\xi''(E')$ for $E' \in \NN(E)_k$ (resp. $E' \in \NN(E)_{>k}$).
Let \[\Xiij i j :=\Xi''(E_1)_i \cap \Xi''(E_2)_j.\] We claim that we
may express the rescaling factors $\Xi''(E_1)_{>0} \cup
\Xi''(E_2)_{>0}$ in terms of the other $N_{r-2}+2$ rescaling factors
$\{\eta''_1,\eta''_2\} \cup \Xiij 0 0$.

We will prove this for $r\in \{6,7\}$ as follows: The
$2\cdot(N_r-N_{r-1}-1)$ equations $g_{D,i}$ are homogeneous of degree
$D$ with respect to the $\Pic(S_r)$-grading of the variables
$\xi''(E)$, and we are interested only in the solutions where all
$\xi''(E)$ are non-zero. Because of the action of $\Tns(S_r)$ on the
rescaling factors and as $E_1, \dots, E_r$ are part of a basis of
$\Pic(S_r)$, we may assume $\eta''_1=\dots=\eta''_r = 1$.

Consider a $(1)$-ruling $D=E_1+E$ such that $(E_2,E)=0$. Then \[D =
E'_1+E''_1 = \dots = E'_{r-3}+E''_{r-3} = E_1+E = E_2+E'\] are the
$r-1$ possibilities to write $D$ as the sum of two intersecting
$(-1)$-curves. Here, $E'_i, E''_i \in \NN(E_1)_0 \cap \NN(E_2)_0$. As
in the discussion before Lemma~\ref{lem:equations}, we may set up the
equations $F_{D,j}$ such that $g_{D,1}$ and $g_{D,2}$ express
$\xi''(E)$ and $\xi''(E')$ directly as a linear combination of
$\xi''(E'_i)\xi''(E''_i)$. This expresses all $\xi(E) \in \Xiij 1 0$
and all $\xi''(E') \in \Xiij 0 1$ in terms of variables in $\Xiij 0 0$.

For a $(1)$-ruling $D=E_1+E$ such that $(E_2,E)=1$, we have
\[D=E'_1+E''_1=\dots=E'_{r-2}+E''_{r-2} = E_1+E,\] where we may assume
$(E_2,E'_i)=0$ and $(E_2,E''_i)=1$. Since $(E_1,E'_i)=(E_1,E''_i)=0$,
we have $\xi''(E'_i) \in \Xiij 0 0$ and $\xi''(E''_i) \in \Xiij 0 1$.
Using the previous findings to express $\xi''(E''_i)$ in terms of
variables $\Xiij 0 0$, the equation $g_{D,1}$ results in a condition
on the variables $\Xiij 0 0$, while $g_{D,2}$ expresses $\xi''(E) \in
\Xiij 1 1$ in terms of these variables.

In the case $r=7$, the equation $g_{D,2}$ for $D=E_1+C_2$ expresses
$\lambda''_2 \in \Xiij 1 2$ in terms of variables in $\Xiij 0 1$, and
$g_{D,1}$ gives a further condition on these variables. Furthermore,
$g_{D,2}$ for the $(1)$-ruling $D=E_2+C_1$ expresses $\lambda''_1 \in
\Xiij 2 1$ in terms of $\Xiij 1 0$, while $g_{D,1}$ gives a further
condition on them. Substituting the expressions for $\Xiij 0 1$
respectively $\Xiij 1 0$ in terms of $\Xiij 0 0$, we get expressions
for $\lambda''_2$ and $\lambda''_1$, while the $g_{D,1}$ result in
further condition on $\Xiij 0 0$.

We summarize this in the following lemma. For its proof, it remains to
show in the following sections that the expressions for
$\Xi''(E_1)_{>0} \cup \Xi''(E_2)_{>0}$ are non-zero, and that the
further conditions vanish.

\begin{lemma}\label{lem:factors}
  We can write the $N_r-N_{r-2}-2$ rescaling factors in the set
  $\Xi''(E_1)_{>0} \cup \Xi''(E_2)_{>0}$ as non-zero expressions in
  terms of $N_{r-2}+2$ rescaling factors $\Xiij 0 0 \cup \{\xi''(E_1),
  \xi''(E_2)\}$. With this, the $N_r-2N_{r-1}+N_{r-2}$ further
  conditions on the rescaling factors are trivial.
\end{lemma}

For an open subset of the $N_{r-2}+2$ parameters $\{\eta''_1,
\eta''_2\} \cup \Xiij 0 0$, all rescaling factors are non-zero, so we
obtain good rescalings, which proves Theorem~\ref{thm:main} once the
proof of Lemma~\ref{lem:factors} is completed.

\section{Degree 3}\label{sec:degree_3}

In this section, we prove Lemma~\ref{lem:factors} for $r=6$ by solving
the system of equations on the rescaling factors of
Lemma~\ref{lem:equations}: For each $(1)$-ruling $D \in \MM_6$, we
determine the coefficients of the equation $p_D$ defining $H_6$ as in
\eqref{eq:defining_Hr}, and find the coefficients $\alpha_j, \beta_j$
of $F_{D,j}$ defining $\AA(S_6)$ as in \eqref{eq:defining_FDj} in the
list in \cite[Section~3]{math.AG/0603111}. This allows us to write
down the 20 equations $g_{D,i}$ on the rescaling factors $\xi''(E)$
explicitly.  Let
\[\gamma_1 := ad-bc, \qquad \gamma_2:=(a-1)(d-1)-(b-1)(c-1)\] for
simplicity.
\begin{align*}
  g_{E_1+m_{1,2},1} &= -\eta''_3 \mu''_{2,3} - \eta''_4 \mu''_{2,4}
  - b \eta''_5 \mu''_{2,5} - d \eta''_6 \mu''_{2,6},\\
  g_{E_1+m_{1,2},2} &= \eta''_1 \mu''_{1,2} + \eta''_4 \mu''_{2,4}
  + \eta''_5 \mu''_{2,5} + \eta''_6 \mu''_{2,6},\displaybreak[0]\\
  g_{E_1+m_{1,3},1} &= -\eta''_1 \mu''_{1,3} + \eta''_4 \mu''_{3,4}
  + \eta''_5 \mu''_{3,5} + \eta''_6 \mu''_{3,6},\\
  g_{E_1+m_{1,3},2} &= \eta''_2 \mu''_{2,3} + \eta''_4 \mu''_{3,4}
  + a \eta''_5 \mu''_{3,5} + c \eta''_6 \mu''_{3,6},\displaybreak[0]\\
  g_{E_1+m_{1,4},1} &= -\eta''_1 \mu''_{1,4} + \eta''_3 \mu''_{3,4}
  + (b - 1) \eta''_5 \mu''_{4,5} + (1-d) \eta''_6 \mu''_{4,6},\\
  g_{E_1+m_{1,4},2} &= -\eta''_2 \mu''_{2,4} + \eta''_3 \mu''_{3,4}
  + (b-a) \eta''_5 \mu''_{4,5} + (c - d) \eta''_6 \mu''_{4,6},\displaybreak[0]\\
  g_{E_1+m_{1,5},1} &= -\eta''_2 \mu''_{2,5} + a/b \eta''_3 \mu''_{3,5}
  + (a - b)/b \eta''_4 \mu''_{4,5} + \gamma_1/b \eta''_6 \mu''_{5,6},\\
  g_{E_1+m_{1,5},2} &= -\eta''_1 \mu''_{1,5} + 1/b \eta''_3 \mu''_{3,5}
  + (1-b)/b \eta''_4 \mu''_{4,5} + (d-b)/b \eta''_6 \mu''_{5,6},\displaybreak[0]\\
  g_{E_1+m_{1,6},1} &= -\eta''_1 \mu''_{1,6} + 1/d \eta''_3 \mu''_{3,6}
  + (d - 1)/d \eta''_4 \mu''_{4,6} + (b - d)/d \eta''_5 \mu''_{5,6},\\
  g_{E_1+m_{1,6},2} &= -\eta''_2 \mu''_{2,6} + c/d \eta''_3 \mu''_{3,6}
  + (d-c)/d \eta''_4 \mu''_{4,6} -\gamma_1/d \eta''_5 \mu''_{5,6},\displaybreak[0]\\
  g_{E_1+Q_2,1} &=  a(c - d)\eta''_1\lambda''_2 + (d - 1)\eta''_2\lambda''_1 - \mu''_{3,4}\mu''_{5,6} + \mu''_{3,6}\mu''_{4,5},\\
  g_{E_1+Q_2,2} &=  \gamma_1 \eta''_1\lambda''_2 + (b - d)\eta''_2\lambda''_1 - \mu''_{3,5}\mu''_{4,6} - \mu''_{3,6}\mu''_{4,5},\displaybreak[0]\\
  g_{E_1+Q_3,1} &=  b(c - d)\eta''_1\lambda''_3 + (c - 1)\eta''_3\lambda''_1 - \mu''_{2,4}\mu''_{5,6} + \mu''_{2,6}\mu''_{4,5},\\
  g_{E_1+Q_3,2} &=  \gamma_1 \eta''_1\lambda''_3 + (a - c)\eta''_3\lambda''_1 - \mu''_{2,5}\mu''_{4,6} - \mu''_{2,6}\mu''_{4,5},\displaybreak[0]\\
  g_{E_1+Q_4,1} &=  bc\eta''_1\lambda''_4 + (bc - b - c + 1)\eta''_4\lambda''_1 - \mu''_{2,3}\mu''_{5,6} + \mu''_{2,6}\mu''_{3,5},\\
  g_{E_1+Q_4,2} &=  (ad - bc)\eta''_1\lambda''_4 + \gamma_2 \eta''_4\lambda''_1 + \mu''_{2,5}\mu''_{3,6} - \mu''_{2,6}\mu''_{3,5},\displaybreak[0]\\
  g_{E_1+Q_5,1} &=  (d-c)\eta''_1\lambda''_5 + \gamma_2\eta''_5\lambda''_1 - \mu''_{2,4}\mu''_{3,6} + \mu''_{2,6}\mu''_{3,4},\\
  g_{E_1+Q_5,2} &=  c\eta''_1\lambda''_5 + (a-c)(1-b)\eta''_5\lambda''_1 - \mu''_{2,3}\mu''_{4,6} - \mu''_{2,6}\mu''_{3,4},\displaybreak[0]\\
  g_{E_1+Q_6,1} &= (b-a) \eta''_1 \lambda''_6 + \gamma_2 \eta''_6
  \lambda''_1
  + \mu''_{2,4} \mu''_{3,5} - \mu''_{2,5} \mu''_{3,4},\\
  g_{E_1+Q_6,2} &= a \eta''_1 \lambda''_6 + (c-a)(d-1)
  \eta''_6 \lambda''_1
  - \mu''_{2,3} \mu''_{4,5} + \mu''_{2,5} \mu''_{3,4}.
\end{align*}
As explained in the previous section, we may assume
$\eta_1''=\dots=\eta_6''=1$.

Recall the discussion before the definition \eqref{eq:defining_FDj} of
$F_{D,j}$ and before Lemma \ref{lem:equations}. If we had chosen
$\xi_4 = \eta_i\lambda_1$ and $\xi_5 = \eta_1\lambda_i$ when writing
down the equations $F_{E_1+Q_i,j}$ in \cite[Section~3]{math.AG/0603111}, then the resulting $g_{E_1+Q_i,2}$ would give
$\lambda''_i$ directly as a quadratic expression in terms of
$\mu''_{j,k}$. Furthermore, each of the five $g_{E_1+Q_i,1}$ would
express $\lambda''_1$ as a quadratic equation in $\mu''_{j,k}$. Of
course, we get the same result by solving the equivalent system of
equations $g_{E_1+Q_i,j}$ as listed above.

The equations $g_{E_1+m_{1,i},j}$ for $i \in \{3, \dots, 6\}$ and
$g_{E_1+Q_2,j}$ allow us to express the variables $\mu''_{1,i}$,
$\lambda''_2$ in $\Xiij 1 0$ and $\mu''_{2,i}$, $\lambda''_1$ in
$\Xiij 0 1$ in terms of the six variables $\mu''_{3,4}, \dots,
\mu''_{5,6} \in \Xiij 0 0$. (As we have set the remaining elements
$\eta''_3, \dots, \eta''_6$ of $\Xiij 0 0$ to the value $1$, they do
not occur in these expressions.)

With $\gamma_3:=d(a-c)(1-b)-c(b-d)(1-a)$, we obtain:
\begin{align*}
  \mu''_{1,3} =& \mu''_{3,4} + \mu''_{3,5} + \mu''_{3,6},\\
  \mu''_{2,3} =& - \mu''_{3,4} - a \mu''_{3,5} - c \mu''_{3,6},\displaybreak[0]\\
  \mu''_{1,4} =& \mu''_{3,4} + (b-1) \mu''_{4,5} + (1-d) \mu''_{4,6},\\
  \mu''_{2,4} =& \mu''_{3,4} + (b - a) \mu''_{4,5} + (c-d) \mu''_{4,6},\displaybreak[0]\\
  \mu''_{1,5} =& 1/b \mu''_{3,5} + (1 - b)/b \mu''_{4,5} + (d - b)/b \mu''_{5,6},\\
  \mu''_{2,5} =& a/b \mu''_{3,5} + (a-b)/b \mu''_{4,5} +\gamma_1/b \mu''_{5,6},\displaybreak[0]\\
  \mu''_{1,6} =& 1/d \mu''_{3,6} + (d-1)/d \mu''_{4,6} + (b-d)/d \mu''_{5,6},\\
  \mu''_{2,6} =& c/d \mu''_{3,6} + (d - c)/d \mu''_{4,6} - \gamma_1/d \mu''_{5,6},\displaybreak[0]\\
  \lambda''_1 =&-\gamma_1/\gamma_3\mu''_{3,4}\mu''_{5,6} -
  a(d-c)/\gamma_3\mu''_{3,5}\mu''_{4,6} - c(b-a)/\gamma_3\mu''_{3,6}\mu''_{4,5},\\
  \lambda''_2 =&(b - d)/\gamma_3\mu''_{3,4}\mu''_{5,6} +
  (1-d)/\gamma_3\mu''_{3,5}\mu''_{4,6} +
  (1-b)/\gamma_3\mu''_{3,6}\mu''_{4,5}.
\end{align*}

For $E \in \{m_{1,2}, Q_3, \dots, Q_6\}$, we consider the remaining
equations $g_{E_1+E,i}$. We can use $g_{E_1+E,2}$ and substitution of
our previous results in order to express $\xi''(E) \in \Xiij 1 1$ in
terms of $\Xiij 0 0$:
\begin{align*}
  \mu''_{1,2} =& \mu''_{3,4} - a/b \mu''_{3,5} - c/d \mu''_{3,6} +
  (a-b)(b-1)/b \mu''_{4,5}\\ &+
  (d-c)(d-1)/d \mu''_{4,6} +(b-d)\gamma_1/(b d) \mu''_{5,6},\\
  \lambda''_3 =&(a - c)/\gamma_3\mu''_{3,4}\mu''_{5,6} + a(1-c)/(b\gamma_3)\mu''_{3,5}\mu''_{4,6} + c(1-a)/(d\gamma_3)\mu''_{3,6}\mu''_{4,5}\\&+ 1/(bd)\mu''_{4,5}\mu''_{4,6} - 1/d\mu''_{4,5}\mu''_{5,6} + 1/b\mu''_{4,6}\mu''_{5,6},\\
  \lambda''_4 =& \gamma_2/\gamma_3\mu''_{3,4}\mu''_{5,6} -
  1/(bd)\mu''_{3,5}\mu''_{3,6} +
  (1-d)(c-d)(a-1)/(d\gamma_3)\mu''_{3,5}\mu''_{4,6}\\ &-
  1/d\mu''_{3,5}\mu''_{5,6} + (1-c)(b-1)(a-b)/(b\gamma_3)\mu''_{3,6}\mu''_{4,5} - 1/b\mu''_{3,6}\mu''_{5,6},\\
  \lambda''_5 = &1/d\mu''_{3,4}\mu''_{3,6} - 1/d\mu''_{3,4}\mu''_{4,6}
  + (b-d)(1-a)\gamma_1/(d\gamma_3)\mu''_{3,4}\mu''_{5,6}\\ &+
  a\gamma_2/\gamma_3\mu''_{3,5}\mu''_{4,6} +(b-1)(a-b)(a-c)/\gamma_3\mu''_{3,6}\mu''_{4,5} -  \mu''_{3,6}\mu''_{4,6},\\
  \lambda''_6 = &-1/b\mu''_{3,4}\mu''_{3,5} -
  1/b\mu''_{3,4}\mu''_{4,5}
  +(1-c)(b-d)\gamma_1/(b\gamma_3)\mu''_{3,4}\mu''_{5,6} \\
  &-\mu''_{3,5}\mu''_{4,5}+(d-1)(c-d)(a-c)/\gamma_3\mu''_{3,5}\mu''_{4,6}
  +c\gamma_2/\gamma_3\mu''_{3,6}\mu''_{4,5}
\end{align*}
Finally, we check by substituting and expanding that the five further
conditions $g_{E_1+E,1}$ are trivial.

Using the restrictions on $a$, $b$, $c$, $d$ imposed by the fact that
$p_1, \dots, p_6$ are in general position (e.g., $a$ must be different
from $b$ and $c$, and all are neither 0 nor 1), we see that
$\mu''_{1,2}, \dots, \mu''_{2,6}, \lambda''_1, \dots, \lambda''_6$ are
non-zero polynomials in $\mu''_{3,4}, \dots, \mu''_{5,6}$.  Therefore,
for an open subset of the $N_4+2=12$ parameters $\eta''_1, \dots,
\eta''_6, \mu''_{3,4}, \dots, \mu''_{5,6}$, all rescaling factors are
non-zero, resulting in good rescalings.

\section{Degree 2}\label{sec:degree_2}

For the proof of Lemma~\ref{lem:factors} for $r=7$, we proceed as in
the case $r=6$ and assume $\eta''_1 = \dots = \eta''_7 = 1$.

Let $D := D^{(1)}_i \in \MM_7$ for $i \in \{3, \dots, 7\}$. We can
arrange $F_{D,1}, \dots F_{D,4}$ in such a way that $g_{D,1}$ and
$g_{D,2}$ express $\mu''_{1,i}$ and $\mu''_{2,i}$ in terms of
$\mu''_{i,j}$ for $j \in \{3, \dots, 7\} \setminus \{i\}$ (see the
discussion before Lemma~\ref{lem:equations}). Similarly, for $i \in
\{3, \dots, 7\}$ and $D := D^{(2)}_{1,2,i} \in \MM_7$, we can arrange
$g_{D,1}$ and $g_{D,2}$ such that they express $\nu''_{1,i}$ and
$\nu''_{2,i}$ linearly in $\nu''_{1,2}$ and of degree 2 in
$\mu''_{3,4}, \dots, \mu''_{6,7}$. This expresses all variables in
$\Xiij 0 1 \cup \Xiij 1 0$ in terms of $\Xiij 0 0$.

Substituting this into an appropriately arranged $g_{D,2}$ for $D :=
D^{(1)}_2 \in \MM_7$ gives $\mu''_{1,2} \in \Xiij 1 1$ in terms of
$\mu''_{3,4}, \dots, \mu''_{6,7}$, and we check that $g_{D,1}$ becomes
trivial.

Now, let $D := D^{(2)}_{1,i,j} \in \MM_7$ for $i<j \in \{3, \dots,
7\}$. We arrange $g_{D,1}$ and $g_{D,2}$ such that they express
$\nu''_{1,i}$ respectively $\nu''_{i,j}$ in terms of $\nu''_{1,j}$ and
expressions of degree 2 in $\mu''_{k,l}$. Using our previous findings,
the first expression turns out trivial, and the second one gives
$\nu''_{i,j} \in \Xiij 1 1$ in terms of $\nu''_{1,2}, \mu''_{3,4}, \dots,
\mu''_{6,7} \in \Xiij 0 0$.

Finally, let $D \in \{D^{(3)}_{1,2}, \dots, D^{(3)}_{1,7},
D^{(3)}_{2,1}\} \subset \MM_7$. We arrange $g_{D,1}$ and $g_{D,2}$
such that the first one is an expression in $\mu''_{j,k}$ and
$\nu''_{j,k}$ which becomes trivial. The second one expresses
$\lambda''_i$ in terms of $\mu''_{j,k}$ and $\nu''_{j,k}$, and we
substitute again our previous results.

This completes the proof of Lemma~\ref{lem:factors} and thus
Theorem~\ref{thm:main}. In total, we obtain good rescalings for an
open subset of a system $N_5+2=18$ parameters \[\eta''_1, \dots,
\eta''_7, \quad \mu''_{3,4}, \dots, \mu''_{6,7}, \quad \nu''_{1,2}.\]
Since it is straightforward to determine the exact expressions for the
remaining 38 rescaling factors in terms of these parameters, and
since the expressions are rather long, we choose not to list them
here.

\begin{remark}
  In principle, it would be possible to consider the conditions
  $g_{D,i}$ for all $(1)$-rulings without reducing to the subset
  $\MM_r$ as we did in Section~\ref{sec:homogeneous_spaces}. While
  this is doable in degree 3 with some software help (\texttt{Magma}),
  especially the expressions corresponding to the $(1)$-ruling
  $D^{(5)}_i$ in degree 2 seem to be out of reach for direct
  computations. Furthermore, we would have to embed the relations
  $v^i_i$ corresponding to the $(2)$-ruling $-K_7$, which causes
  further complications.
\end{remark}

\begin{remark}
  For $r \in \{5,6,7\}$, there is a $(N_{r-2}-r+1)$-parameter family
  of images of $\AA(S_r)$ under good embeddings in $H_r$ by
  Remark~\ref{rem:parameters}. The dimension of $\AA(S_r)$ is $r+3$,
  and there is a $(2\cdot(r-4))$-parameter family of smooth Del Pezzo
  surfaces $S_r$ of degree $9-r$. The dimension of $H_r$ is
  $N_{r-1}+1$.

  In fact, \cite[Section~6.3]{popov_diplom} shows that the closure of
  the union of all these images for all Del Pezzo surfaces of degree 4
  equals $H_5$.

  For $r=6$, by comparing the dimensions and numbers of parameters, a
  similar result seems possible. However, for $r=7$, we have
  \[(N_{r-2}-r+1)+(r+3)+2\cdot(r-4) = 26,\] while $H_7$ has dimension
  $N_{r-1}+1=28$.  Consequently, the closure of the union of the
  corresponding images, over all Del Pezzo surfaces of degree 2, under
  all good embeddings cannot be $H_7$ for dimension reasons.
\end{remark}

\bibliographystyle{alpha}

\bibliography{cox_embedding}

\end{document}